\documentclass[runningheads,a4paper]{llncs}

\usepackage{amssymb,relsize,amsmath,longtable,dcolumn,mathrsfs,subfigure,bm,epstopdf}
\usepackage[vlined,ruled]{algorithm2e}
\usepackage{graphicx}

\usepackage{url}
\urldef{\mailsa}\path|{jan.baetens, bernard.debaets}@ugent.be|

\def\T{{\cal T}}

\newcolumntype{d}[1]{D{.}{.}{#1}}
\newcolumntype{.}{D{.}{.}{-1}}

\begin{document}

\mainmatter  % start of an individual contribution

% first the title is needed
\title{Introducing Lyapunov profiles of\\ cellular automata}

% a short form should be given in case it is too long for the running head
\titlerunning{Lyapunov profiles of cellular automata}

% the name(s) of the author(s) follow(s) next
%
% NB: Chinese authors should write their first names(s) in front of
% their surnames. This ensures that the names appear correctly in
% the running heads and the author index.
%
\author{Jan M.\ Baetens \inst{1}
\thanks{Corresponding author}%%
\and Janko Gravner \inst{2}}
\authorrunning{Jan M.\ Baetens \and Janko Gravner}
% (feature abused for this document to repeat the title also on left hand pages)

% the affiliations are given next; don't give your e-mail address
% unless you accept that it will be published
\institute{KERMIT, Department of Mathematical Modelling, Statistics and Bioinformatics,
Ghent University, Coupure links 653, Gent, Belgium\\ \email{jan.baetens@ugent.be} \and Mathematics Department, University of California, Davis, CA 95616, USA\\ \email{gravner{@}math.ucdavis.edu} 
%\url{http://www.springer.com/lncs}
}
%
% NB: a more complex sample for affiliations and the mapping to the
% corresponding authors can be found in the file "llncs.dem"j
% (search for the string "\mainmatter" where a contribution starts).
% "llncs.dem" accompanies the document class "llncs.cls".
%

%\toctitle{Towards adequate measures grasping CA behavior}
%\tocauthor{Authors' Instructions}
\maketitle

\begin{abstract}
In line with the stability theory of continuous dynamical systems, Lyapunov exponents of cellular automata (CAs) have been conceived two decades ago to quantify to what extent their dynamics changes following a perturbation of their initial configuration. More precisely, Lyapunov exponents of CAs have either been understood as the rate by which the resulting defect cone widens as these dynamical systems are evolved from a perturbed initial configuration, or as the rate by which defects accumulate during their evolution. The former viewpoint yields insight into the extent of the affected region, whereas the latter tells us something about the intensity of the defect propagation. In this paper, we will show how these viewpoints can be united by relying on Lyapunov profiles of CAs.
\end{abstract}

%%% Please eliminate the following two lines from your file.
%\clearpage
%\tableofcontents

\section{Introduction}
\label{sec:in}
Ever since their postulation by von Neumann \cite{neumann1948}, researchers have often been struck by the intriguing spatio-temporal dynamics that cellular automata (CAs) can bring forth in spite of their intrinsically simple nature \cite{grassberger1984,wolfram1984b}. In order to arrive at a comprehensive analysis of their behavior, different measures have been proposed, such as Lyapunov exponents \cite{bagnoli1992,shereshevsky1991} and others \cite{ilachinski2001}. With respect to Lyapunov exponents, two distinct viewpoints have emerged during the last two decades. The first one traces back to a suggestion by Wolfram \cite{wolfram1984b}, which was formalized by Shereshevsky \cite{shereshevsky1991}, and later on adopted by several authors \cite{courbage2006,tisseur2000}, and which identifies Lyapunov exponents of elementary CAs with the rates by which the damage front moves to the right and to the left. Consequently, every elementary CA is characterized by means of two exponents, which obviously might differ. On the other hand, the viewpoint put forward by Bagnoli et al. \cite{bagnoli1992} involves only one Lyapunov exponent for every CA that quantifies the exponential rate by which defects accumulate if the CA is evolved from a initial configuration with a single defect without accounting for the direction in which defects propagate \cite{bagnoli1992}. 

The definition of Lyapunov exponents by Bagnoli et al. \cite{bagnoli1992}  actually constitutes the counterpart of the one by Shereshevsky \cite{shereshevsky1991} because accumulation is neglected while directionality is considered in the latter, whereas the opposite holds for Lyapunov exponents according to Bagnoli at al. In other words, the former viewpoint gives insight into the intensity of the defect propagation, whereas the latter indicates how fast an initial perturbation spreads as the CA evolves and therefore how fast the defect cone widens as the CA evolves. Clearly, each viewpoint has its advantages, but uniting them in a single concept  enables a comprehensive stability analysis of CAs. 

In this paper, we will show how so-called Lyapunov profiles of CAs constitute such a uniting framework that allows for assessing both the intensity of the damage propagation and the growth rate of the defect cone. 

\section{Lyapunov profiles of cellular automata}
\subsection{Preliminaries}
In the framework of this paper, we will restrict ourselves to elementary CAs which can be conveniently represented by a quintuple $\mathscr{C}=\left\langle \T,S,s,N,\Omega\right\rangle$. Therein,  $\mathcal{T}$ denotes a countably infinite tessellation of a  $1$-dimensional Euclidean space, consisting of consecutive intervals $c_i$,  $i\in\mathbb{N}$, typically referred to as cells, and $S$ is a set of two states, here $S=\left\{0,1\right\}$. Further,  the output function $s:\,\T\times\mathbb{N}\to S$ gives the state value of cell $c_i$ at the $t$-th discrete time step, $N(c_i)$ is the ordered list of neighbors of $c_i$, \emph{i.e.\ }$N(c_i)=\left(c_{i-1},c_i,c_{i+1}\right)$, and finally the transition function $\Omega:\mathbb{N}\to S$ that governs the dynamics of each cell $c_i$ reads
$$
s(c_i,t+1)=\Omega\left(s(c_{i-1},t),s(c_i,t),s(c_{i+1},t)\right)\,.
$$

\subsection{Rationale}
The Lyapunov exponent of a CA according to  Bagnoli et al. \cite{bagnoli1992} is computed by tracking all defects that emerge from the introduction of a single defect in the CA's initial configuration, \emph{i.e.\ }by flipping one of the cells' state, as follows:
\begin{equation}
\lambda=\displaystyle\lim\limits_{t\rightarrow\infty}\dfrac{1}{t}\log\left(\dfrac{\epsilon_t}{\epsilon_0}\right)\,,\label{lya}
\end{equation}
where $\epsilon_t$ denotes the total number of defects at the $t$-th step and a defect is defined as the smallest possible perturbation of the CA, such that it may be conceived as an object that flips that state of the cell in which it resides. Upon introducing such a defect to one of the CA's cells, it will propagate and affect increasingly more cells because the dynamics of the CA is governed by its neighborhood configurations. 
The quantity $\epsilon_t$ should be computed by evolving both the initial configuration $s_0$ and its perturbed counterpart $s_0^*$, for which it holds that $d(s_0,s_0^*)=1$, for one time step, after which a replica $s^*_i(\cdot,1)$ of $s(\cdot,1)$ is created for every cell $c_i$ for which it holds that $s(c_i,1)\neq s^*(c_i,1)$ in such a way that $d(s(\cdot,1),s^*_i(\cdot,1))=1$. Then, all replicas are evolved one more time step and the resulting configurations $s^*_i(\cdot,2)$ are again compared with $s(\cdot,2)$, such that a new set of replicas can be constructed, after which this procedure is repeated until the CA evolution halts. So, when computing the total number of defects up to a given time step, both the position and multiplicity of defects is known, but this information is lost when computing the Lyapunov according to Eq.~\eqref{lya}. 

Clearly, if we redefine $\epsilon_t$ as a vector $\bm{\epsilon}_t$ whose $i$-th element $\epsilon^i_t$ is the number of defects in cell $c_i$ at a given step, we have a construct that incorporates both the multiplicity and position of defects. For comprehensiveness, it should be mentioned that Bagnoli et al. \cite{bagnoli1992} do incorporate such a damage vector to define the Lyapunov exponent of elementary CAs, but they discard its information content as their further analysis is based upon the sum of its elements. Besides, they call the quantity given by Eq.~\eqref{lya} the maximum Lyapunov exponent (MLE) of a CA, by analogy to continuous $n$-dimensional dynamical systems where the terms refers to the largest exponent in the Lyapunov spectrum that contains the perturbations, and rates of separation, in all $n$ different directions. By contrast, in the case of binary CAs that Bagnoli et al. \cite{bagnoli1992} consider, they define only one exponent, thus the meaning of the term `maximum' in their definition is unclear.

In this paper we interpret a CA on an array of $n$ cells as an $n$-dimensional system, thus it may be anticipated that a CA evolution generates a spectrum of Lyapunov exponents, one per each cell of the array.  The largest of these is naturally called the MLE. As $\bm{\epsilon}_t$  quantifies how the number of defects grows in each possible direction, \emph{i.e.\ }in each cell, and therefore provides insight into the damage propagation resulting from a single defect in each direction,  we proceed to take a closer look at what exactly we might learn from $\bm{\epsilon}_t$ . 

\subsection{Lyapunov profiles}
At this point, we shall define the finite-time Lyapunov profile $P$ of a CA as
$$
\Lambda_T=\dfrac{1}{T}\,\log\left(\bm{\epsilon}_T\right)
$$
where $\log$ is applied element wise. Hence, the Lyapunov profile consists of the time-averaged expansion rates in each of the possible directions. The elements of $\Lambda_T$ may be understood intuitively as the time-averaged exponential rates by which the number of defects grows in the cells of the CA, and we could look at the largest among them as the MLE of the CA. Besides, the width of the defect cone is given by $\left|\left\{i\mid\log\epsilon_T^i\neq-\infty\right\}\right|$. At this point it should also be mentioned that the defect distribution does not translate directly to the distribution of damages in the CA configuration space. In the remainder, we will refer to $\frac{1}{T}\log(\epsilon_T^i)$ as the $i$-th finite-time Lyapunov exponent of a CA, denoted $\lambda_T^i$. 

Given the fact that the propagation of defects is constrained by the size of the neighborhood $N(c_i)$ \cite{bagnoli1992}, there must exist an upper bound on $\Lambda_T$'s elements. Recalling that in the worst-case scenario a defect in a cell $c_i$ will propagate to each of its three neighbors if the CA is evolved for one time step, and this occurs for all cells in $\T$, it is easy to see that the upper bound on $\Lambda_T$ is given by the $T+1$-th row of the trinomial triangle in such a way that 
$$
\epsilon^i_T=\binom{T}{i-i^*_0}_2\,,
$$
where $i^*_0$ denotes the index of the cell that was perturbed at $t=0$.  The sum of its elements on the $T+1$-th row are $3^T$, which is the upper bound on the total number of defects at $T$, which thereby indicates why the MLE according to \cite{bagnoli1992} can be at most $\log(3)$. 

The Lyapunov profiles of rules 105 and 150 at $T$ are completely defined by the $T+1$-th row of the trinomial triangle because all  cells in these CAs pass on defects to all their neighboring cells. Besides, their maximum will be given by
$$
\max\Lambda_T=\dfrac{1}{T}\log\binom{T}{0}_2\,,
$$
or equivalently \cite{andrews1990}:
%\begin{equation}
%\max\Lambda_T=\dfrac{1}{T}\log\left(\sum\limits_{j=0}^{\left\lfloor \frac{T}{2}\right\rfloor}\dfrac{T!}{(j!)^2\,(T-2j)!}\right)\,.
%\end{equation}
%%which can be rewritten as
%%\begin{equation}
%%\max\Lambda_T=\dfrac{1}{T}\log\left(1+\sum\limits_{j=1}^T\exp\left(\log\left(\dfrac{T!}{(j!)^2\,(T-2j)!}\right)\right)\right)\,.
%%\end{equation}
%It can be shown that this equation collapses to:
\begin{equation}
\max\Lambda_T=\dfrac{1}{T}\log\left(1+T!\sum\limits_{j=1}^{\left\lfloor \frac{T}{2}\right\rfloor}\prod\limits_{k=1}^jk^{-2}\prod\limits_{k=1}^{T-2j}k^{-1}\right)\,,
\label{maxlya}
\end{equation}
which provides us with a means to normalize the Lyapunov spectra of elementary CAs, such that a mutual comparison becomes possible.

\section{Results}
\subsection{Experimental setup} 
In the remainder of this paper we will focus our attention on the 88 minimal elementary CAs as defined in \cite{vichniac1990}. For each of them, the propagation of defects emerging from a single defect was tracked for 5000 time steps in a one-dimensional system consisting of 10001 cells, as such mimicking a system of infinite size. In order to enable a meaningful comparison between the Lyapunov profiles of the considered CAs, their elements were normalized by means of the upper bound on the exponential propagation rate as given by Eq.~\eqref{maxlya}. Although the Lyapunov profiles will depend in in many cases on the initial configuration from which they are evolved (Class 2 and 4 rules), we restrict ourselves in this preliminary work to uniformly chosen random initial states and a single initial defect. A way to eliminate some of the dependence is to replace the singleton with a large zone of defects, which makes it very unlikely that defects die by chance.

\subsection{Lyapunov profiles of elementary cellular automata}
\subsubsection{The benchmark rule}
\begin{figure}
  \begin{center}
    %% Don't use epsfile!
		\includegraphics[width=0.3\textwidth]{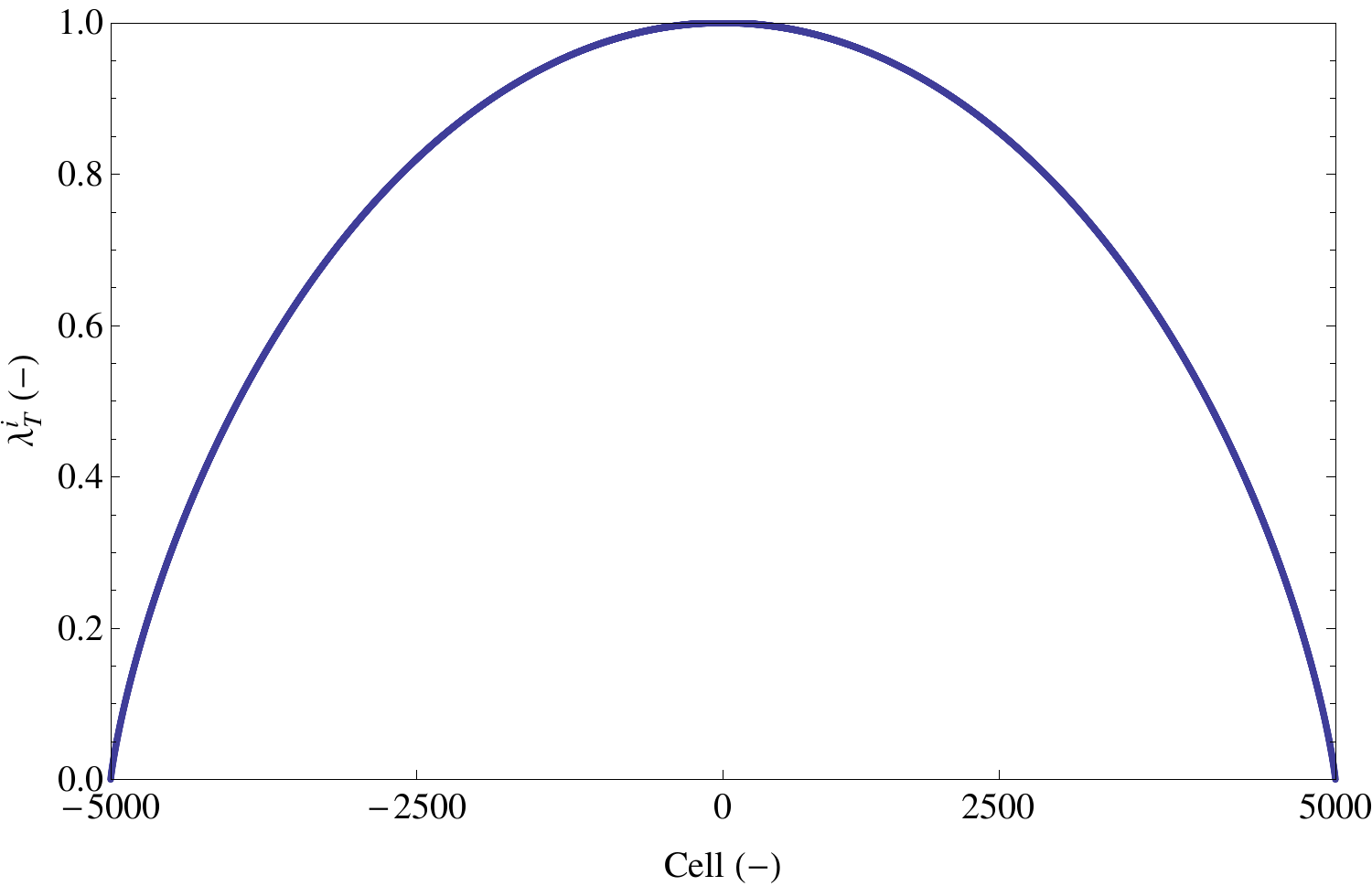}
    \caption{Normalized Lyapunov profile of the ECA rule 150 that was obtained after 5000 time steps for a system of 10001 cells.}
    \label{R150ref}
  \end{center}
\end{figure}

As a kind of benchmark, Fig.~\ref{R150ref} shows the normalized Lyapunov profile of rule 150, which is equivalent to the one of rule 105. For comprehensiveness, the cell indices along the horizontal axis are positioned in such a way that index 0 refers to the cell where initial defect was introduced. As anticipated when deriving the upper bound on the MLE, the normalized profile of this rule reaches a maximum 1 in the initially perturbed cell. The discrepancy between the line $\lambda_T^i=1$ and the profile gives an indication of the delay in the defect propagation that is imposed by the finite speed of information transmission in CAs. Indeed, it takes $t$ time steps before a cell that is $t$ cells apart from the initially perturbed one can start propagating defects. As such, if one would normalize $\Lambda_T$'s elements with respect to the corresponding trinomial coefficient, and not with central trinomial coefficient of the $T+1$-th row of the trinomial triangle, such a normalized Lyapunov profile would result in a horizontal $\lambda_T^i=1$. At this point, we opted not to normalize the Lyapunov profile in this way because one would then somehow lose the information on the delay of the defect propagation that is imposed by the finite speed of information transmission in CAs.

\subsubsection{Class 3 rules}
\begin{figure}
  \begin{center}
    %% Don't use epsfile!
	%	\subfigure[Rule 22]{\includegraphics[width=0.35\textwidth]{R22}}
		\subfigure[Rule 30]{\includegraphics[width=0.45\textwidth]{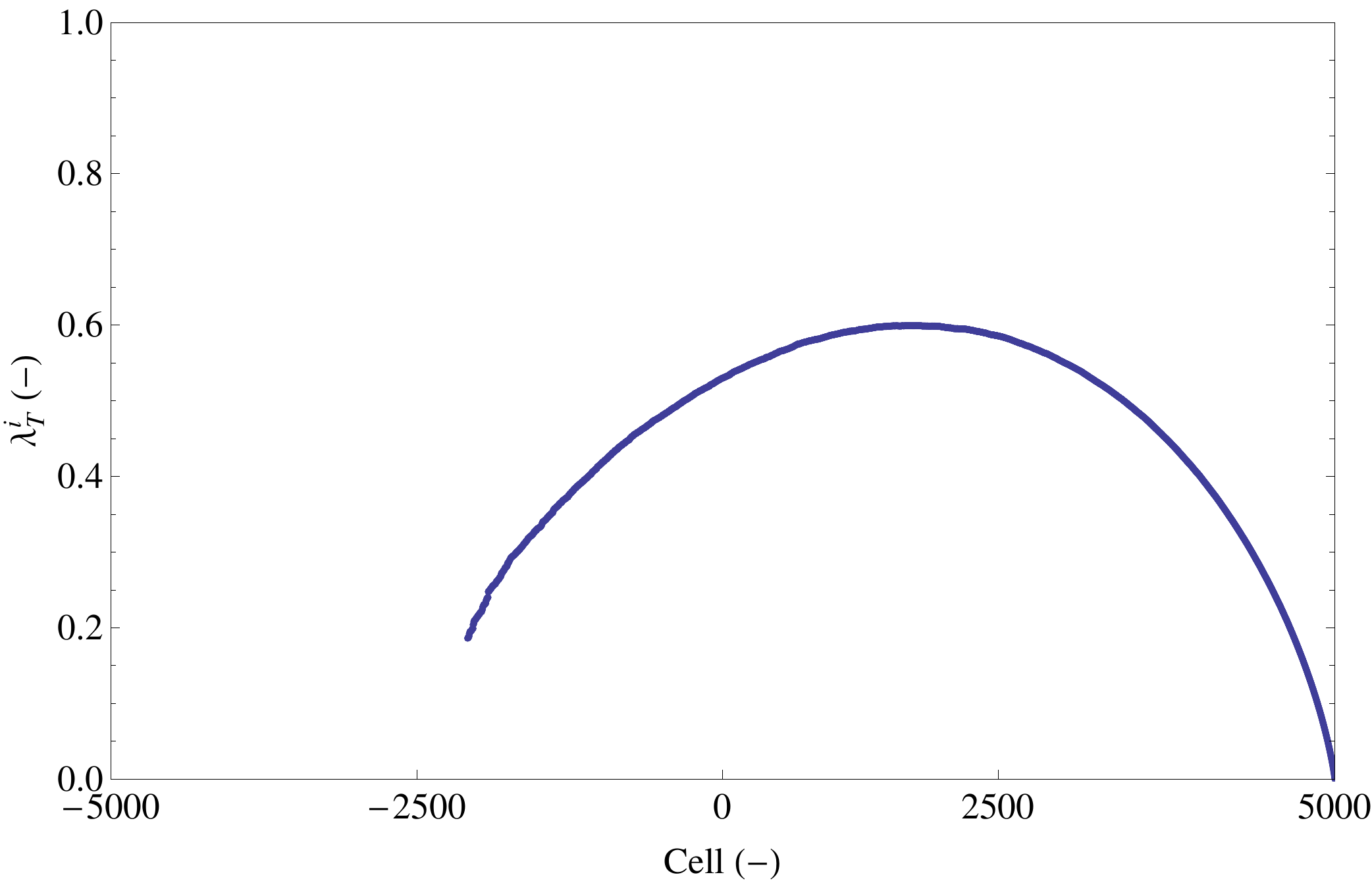}}
				\subfigure[Rule 90]{\includegraphics[width=0.45\textwidth]{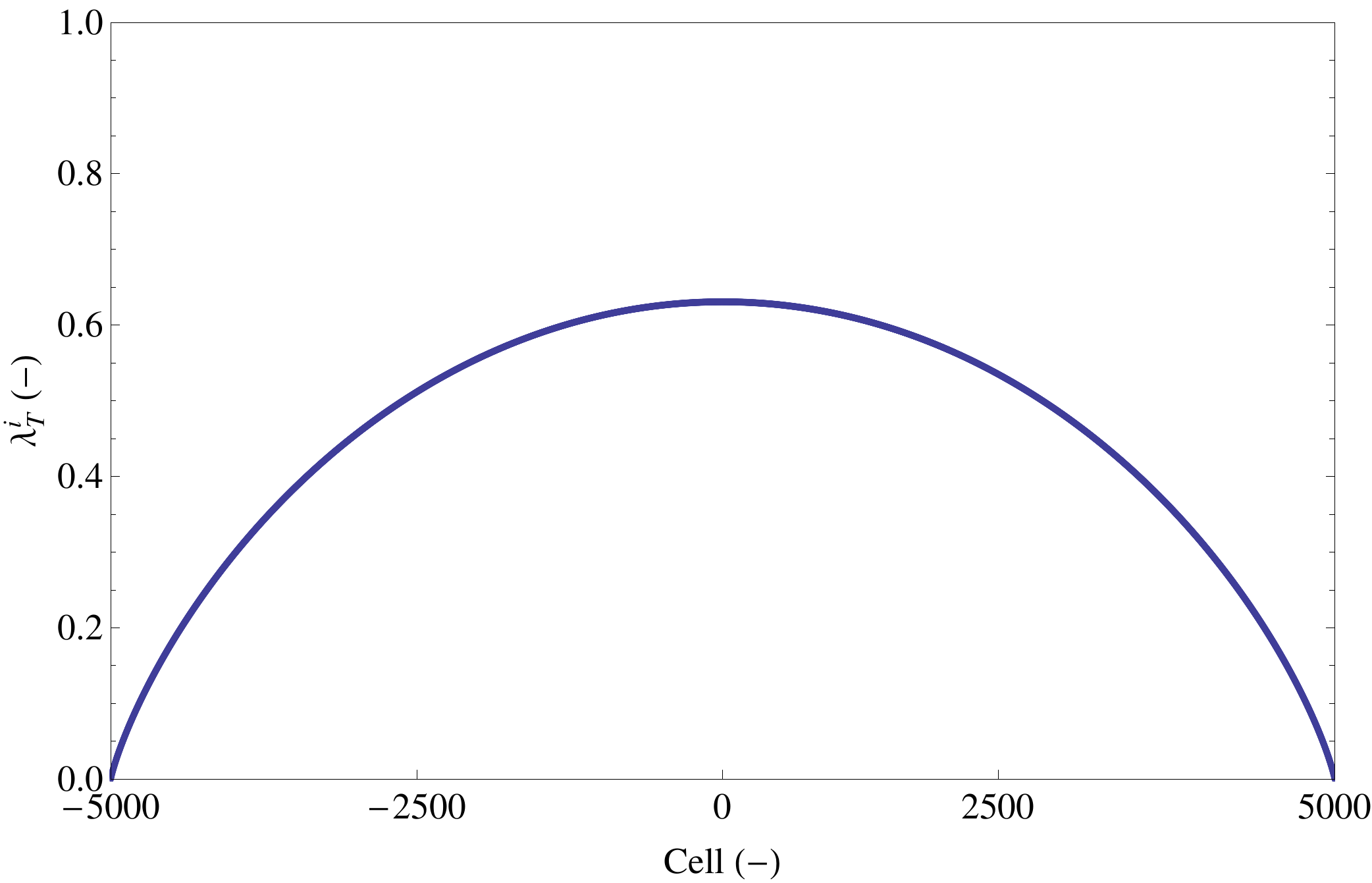}}
		%\subfigure[Rule 122]{\includegraphics[width=0.35\textwidth]{R122}}
    \caption{Normalized Lyapunov profiles of some Class 3 ECAs  that were obtained after 5000 time steps for a system of 10001 cells. }
    \label{Class3}
  \end{center}
\end{figure}

Figure~\ref{Class3} depicts the normalized Lyapunov profiles of some of the most well-studied ECA rules that belong to Wolfram's class 3 \cite{wolfram1984b}. It is clear that the profiles of ECA rules 30 and 90  lie below the one of rule 150, but there are pronounced differences between the profiles of these Class 3 rules. The one of rule 90 shows the most agreement with the one of rule 150, with a very similar shape, but with a maximum of approximately 0.6 instead of 1. It is importance to notice that it has positive values for all cells, which means that the defect cone widens at maximum speed, just as in the case of rule 150, but its intensity of defect propagation is lower. The defect cone for rule 30 is much smaller as it exhibits a sharp transition at $i\approx-2300$ from expansive cells where $\lambda_T^i>0$ to cells that are not yet affected by defects, and hence have $\lambda_T^i=-\infty$. Moreover, the one of rule 30  is  asymmetric, in the sense that its damage cone propagates to the right at maximum speed, whereas it travels to the left at a speed lower than half of this maximum and its maximum is not located at the initially perturbed cell. 

\subsubsection{Class 4 rules}
\begin{figure}
  \begin{center}
    %% Don't use epsfile!
%		\subfigure[Rule 41]{\includegraphics[width=0.35\textwidth]{R41}}
	%	\subfigure[Rule 54]{\includegraphics[width=0.35\textwidth]{R54}}\\
				\subfigure[Rule 106]{\includegraphics[width=0.45\textwidth]{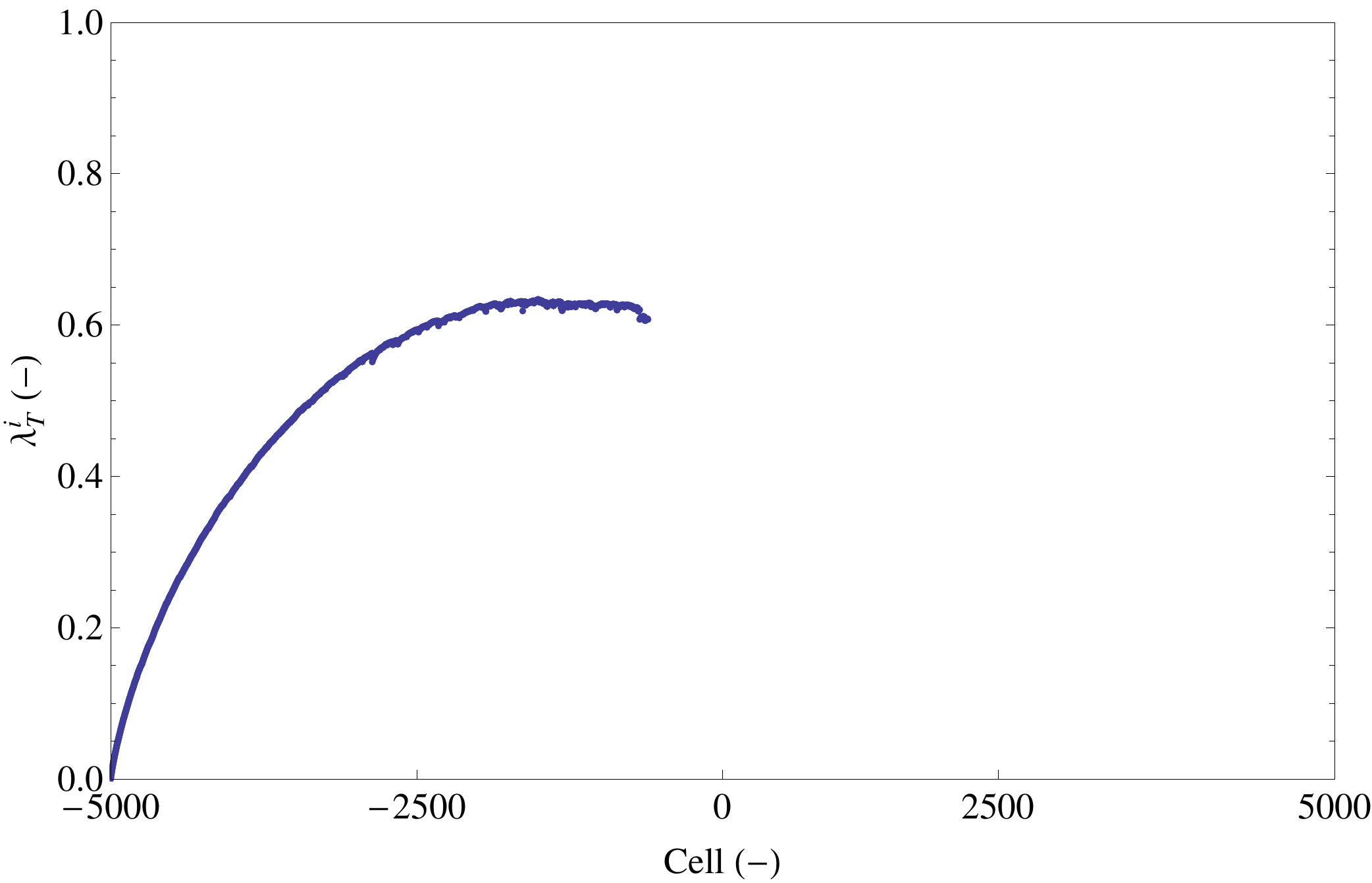}}
		\subfigure[Rule 110]{\includegraphics[width=0.45\textwidth]{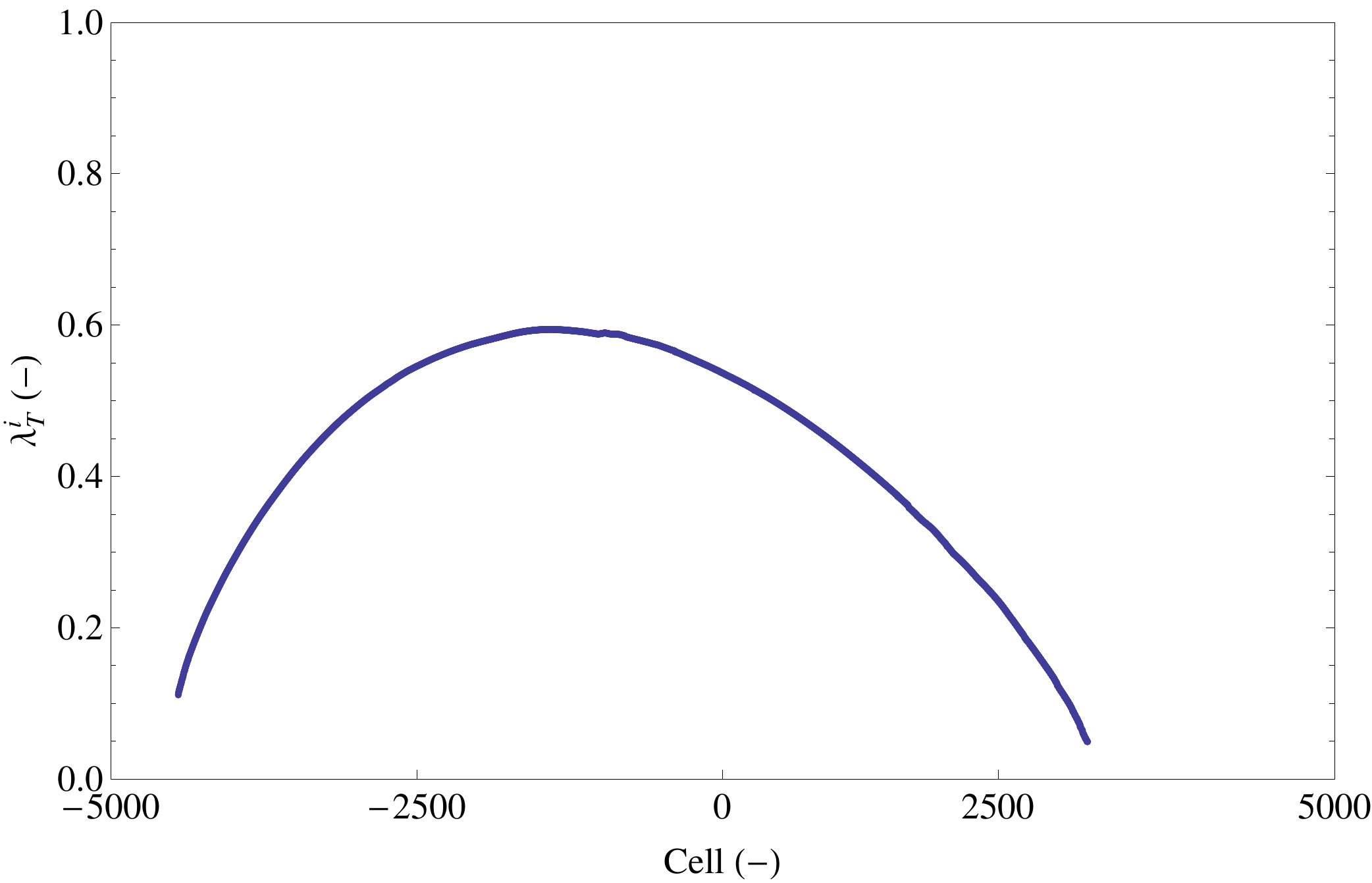}}
    \caption{Normalized Lyapunov profiles of the  Class 4 minimal ECAs  that were obtained after 5000 time steps for a system of 10001 cells. }
    \label{Class4}
  \end{center}
\end{figure}

Figure~\ref{Class4} depicts the normalized Lyapunov profiles of two Wolfram's \cite{wolfram1984b} Class 4 ECA rules that are included among the 88 minimal ECAs. Compared to the Lyapunov profiles of the exemplary Class 3 rules, the ones of the Class 4 rules are more intriguing in the sense that they typically exhibit one or two phase transitions, \emph{i.e.\ }sharp transitions from expansive cells where $\lambda_T^i>0$ to not yet affected cells (rules 106 and 110), are asymmetric (rules 106 and 110), and/or non smooth (rule 106). Besides, just as we observed in the case of rule 30, the maximum of the Lyapunov profile does not necessarily occur at the cell where the defect was introduced at the beginning of the evolution. This skewness is especially pronounced in the case of rule 110. 

\subsubsection{Class 2 rules}
\begin{figure}
  \begin{center}
    %% Don't use epsfile!
		\subfigure[Rule 6]{\includegraphics[width=0.3\textwidth]{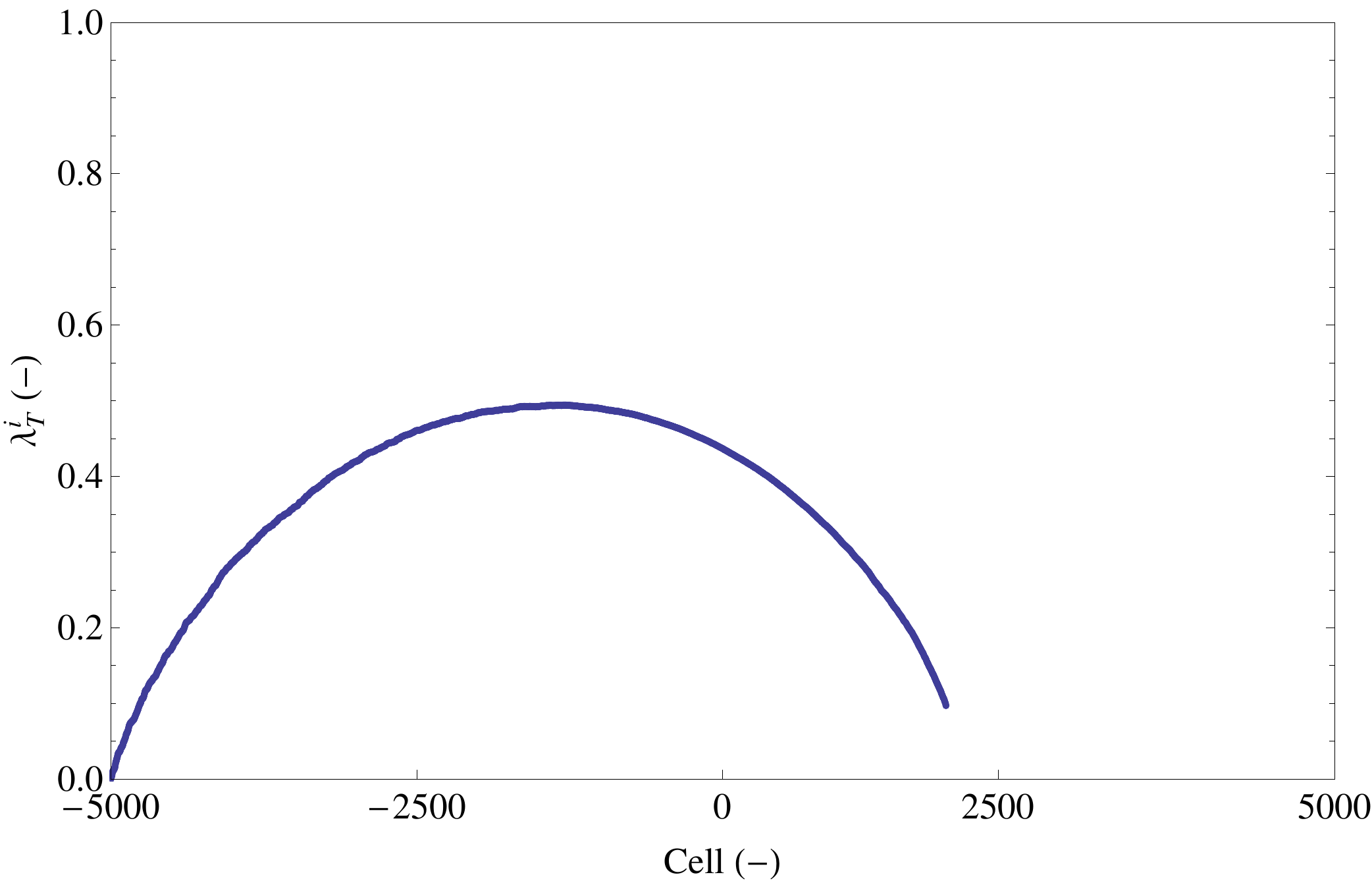}}
		\subfigure[Rule 25]{\includegraphics[width=0.3\textwidth]{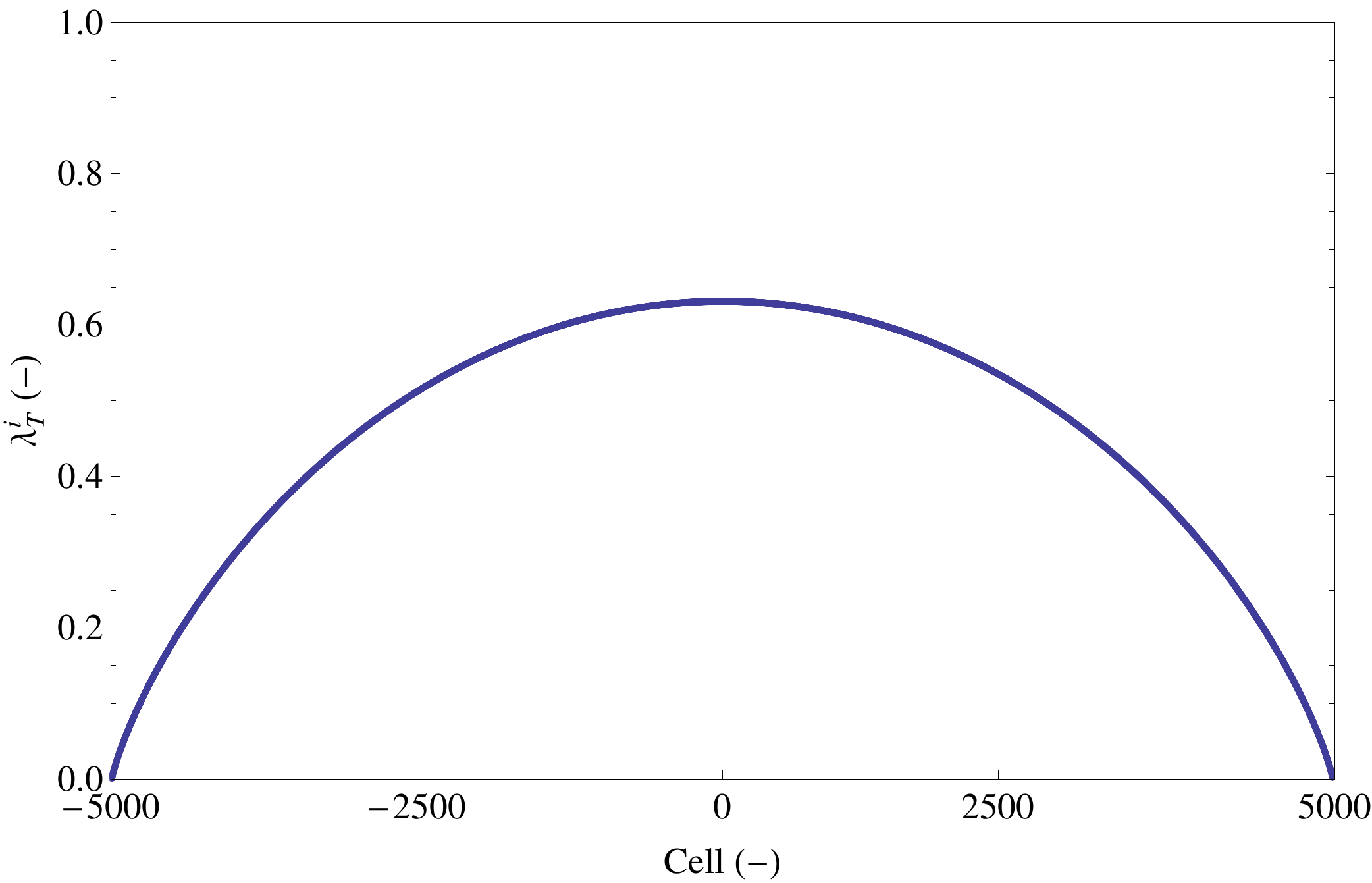}}
				\subfigure[Rule 38]{\includegraphics[width=0.3\textwidth]{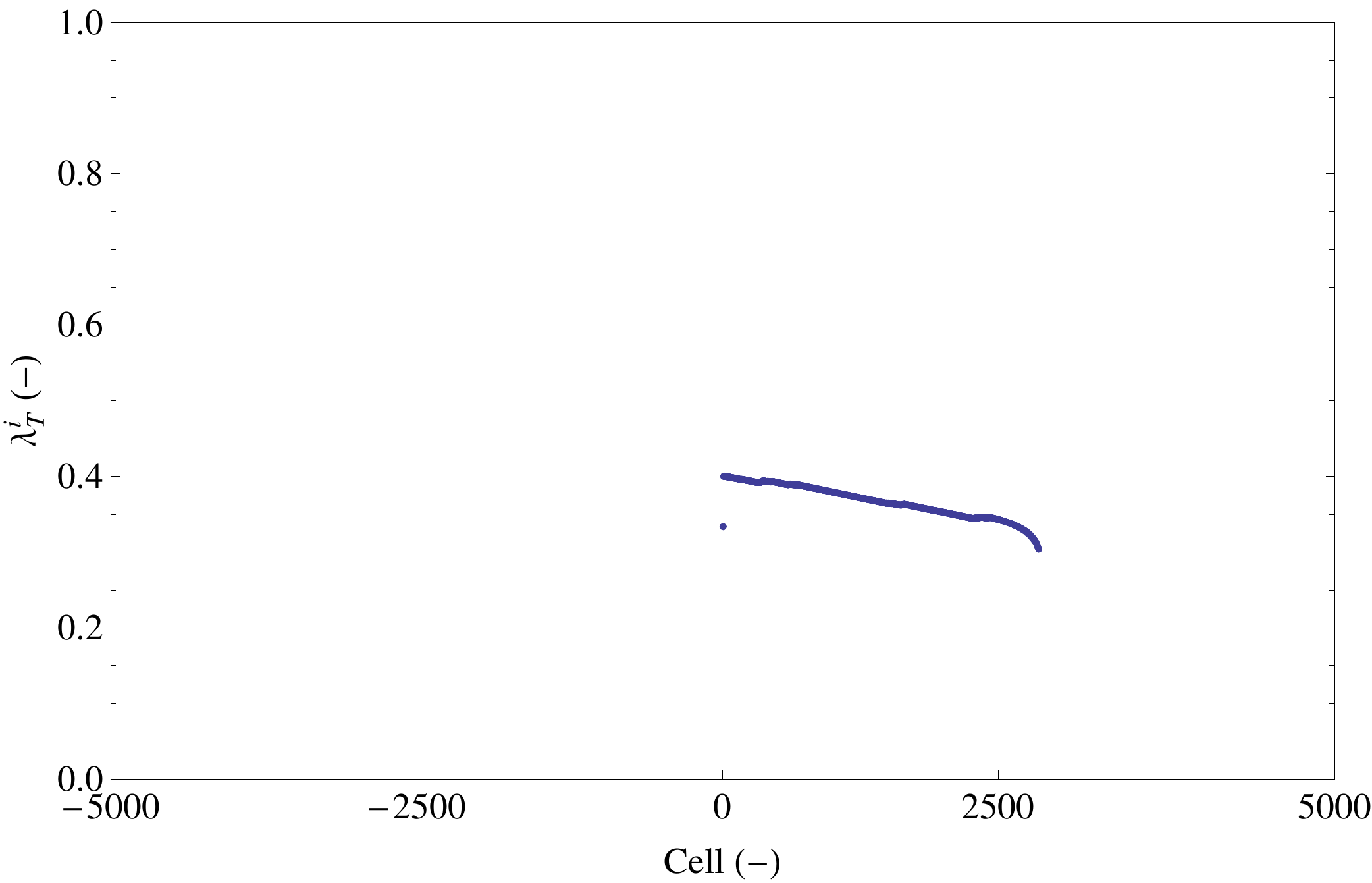}}
		%\subfigure[Rule 57]{\includegraphics[width=0.3\textwidth]{R57}}
						%\subfigure[Rule 62]{\includegraphics[width=0.3\textwidth]{R62}}
		%\subfigure[Rule 134]{\includegraphics[width=0.3\textwidth]{R134}}
    \caption{Normalized Lyapunov profiles of ECA rules 6, 57 and 62 that were obtained after 5000 time steps for a system of 10001 cells. }
    \label{Class2}
  \end{center}
\end{figure}

Although most non-trivial Lyapunov profiles are found for rules belonging to Wolfram's Classes 3 and 4, it should be emphasized that also Class~2 rules can give rise to non-trivial Lyapunov profiles. More precisely, such rules either give rise to  a highly discontinuous profile that are $-\infty$, except in a few cells, or a truly meaningful profile. For instance, within the family of minimal ECAs, rules 28, 33, 37, 73, 108 and 156 give rise to such discontinuous profiles, which indicate that defects can accumulate as the ECAs evolve, but they remain localized and cannot propagate beyond certain limits. We opt not to include exemplary profiles of these rules because they are $-\infty$ almost everywhere. On the other hand, there are a few rules, namely 6, 57 and  62, which give rise to smooth profiles, similar to the ones of the Class 3 rules (Fig.~\ref{Class2}). Other class 2 rules within the studied ECA family lead to a Lyapunov profile that is $-\infty$ everywhere. 

\begin{figure}
  \begin{center}
    %% Don't use epsfile!
\includegraphics[width=0.55\textwidth]{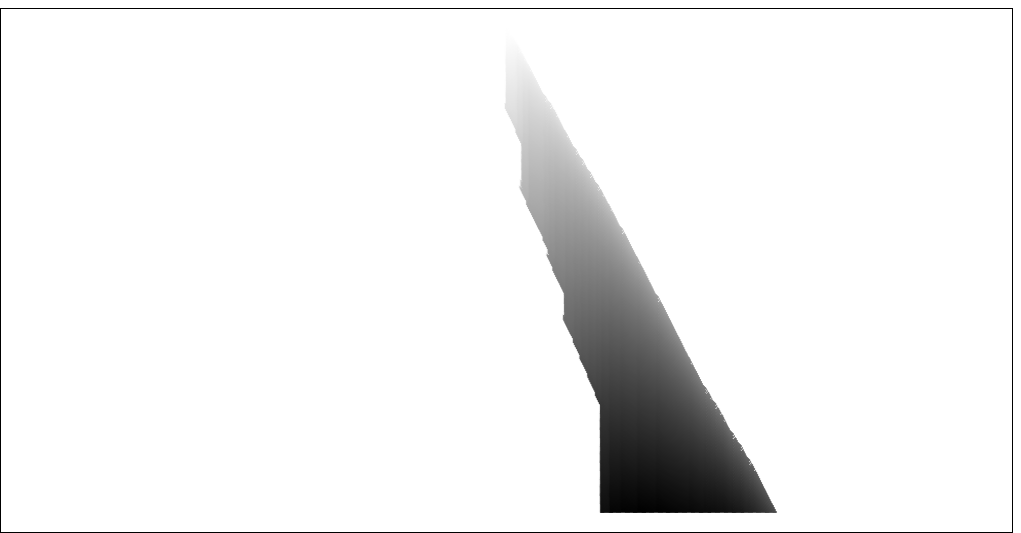}\hspace{0.3cm}\includegraphics[width=0.045\textwidth]{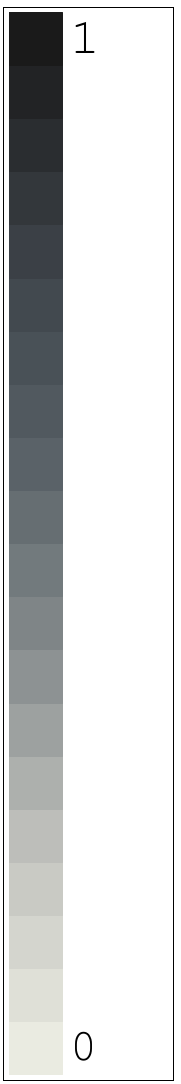}
    \caption{Heatmap of Rule 62 that was obtained after 500 time steps for a system of 1001 cells. }
    \label{Heat}
  \end{center}
\end{figure}

As indicated by the profile of rule 57, even Class 2 rules can give rise to a Lyapunov profile that stretches across the entire light cone, \emph{i.e.\ }$[-T,T]$, but this is the only Class 2 rule that exhibits such a pronounced defect propagation. For the other Class 2 rules that give rise to a non-trivial Lyapunov profile, it is typically lower and narrower than the ones of Class 3 and Class 4 rules. The peculiar Lyapunov profile of rule 62 can be better understood by looking at its corresponding heat map, which displays the time evolution of the normalized Lyapunov profile (Fig.~\ref{Heat}). This clearly shows that the defect propagation is confined to a region of which the right side extends at maximum speed, whereas its left side tends to be pushed inwards substantially from time to time.

\section*{Conclusions and future work}
In this paper it was shown how the distinct viewpoints on Lyapunov exponents of CAs can be united by considering their so-called Lyapunov profiles. These profiles consist of the time-averaged rate of defect accumulation in each of the possible directions that is allowed by the CA's structure. As a next step, we will focus on whether rigorous results can be obtained for the Lyapunov profiles and the therewith defined MLE, their extension to arbitrary two-state families, as well as their formulation in the case of multi-state CAs. 

\section*{Acknowledgements}
This paper has been developed as a result of a mobility stay
funded by the Erasmus Mundus Programme of the European Commission under the
Transatlantic Partnership for Excellence in Engineering – TEE Project. Besides, this work was carried out using the STEVIN Supercomputer Infrastructure at Ghent University, funded by Ghent University, the Flemish Supercomputer Center (VSC), the Hercules Foundation and the Flemish Government – department EWI.

\bibliographystyle{splncs03.bst}
% use the following instead if you encounter problems 
%\bibliographystyle{alpha}
\bibliography{biblioAUTOMATA2014.bib}
\label{sec:biblio}

\end{document}